%%%%%%%   Geometry and Topology Monographs Volume 2 :
%%%%%%%   m2-25.tex  Rudolph for Kirbyschrift

\input gtmacros
\input amsnames
\input amstex
%
          % amstex disables \cal 
\catcode`\@=12        % and leaves @ \active !! caution @ used in math !!
%                                               see hacks above

\input gtmonout
\volumenumber{2}
\volumeyear{1999}
\volumename{Proceedings of the Kirbyfest}
\pagenumbers{555}{562}
\papernumber{25}
\received{31 July 1998}\revised{18 March 1999}
\published{21 November 1999}

%%
%%
%%     Adapt basic layout macros:
%%
\let\\\par
\def\topmatter{\relax}

\let\gttitle\title
\def\title#1\endtitle{\gttitle{#1}}
\let\gtauthor\author
\def\author#1\endauthor{\gtauthor{#1}}
\let\gtaddress\address
\def\address#1\endaddress{\gtaddress{#1}}
\let\gtemail\email
\def\email#1\endemail{\gtemail{#1}}
\def\subjclass#1\endsubjclass{\primaryclass{#1}}
\let\gtkeywords\keywords
\def\keywords#1\endkeywords{\gtkeywords{#1}}
\def\heading#1\endheading{{\def\S##1{\relax}\def\\{\relax\ignorespaces}
    \section{#1}}}
\def\head#1\endhead{\heading#1\endheading}

\def\subhead#1\endsubhead{\sh{#1}}
\def\subsubhead#1\endsubsubhead{\sh{#1}}
\def\specialhead#1\endspecialhead{\sh{#1}}
\def\demo#1{\rk{#1}\ignorespaces}
\def\enddemo{\ppar}
\let\remark\demo
\def\endremark{}

\def\qed{\ifmmode\quad\sq\else\hbox{}\hfill$\sq$\par\goodbreak\rm\fi}  
\def\proclaim#1{\rk{#1}\sl\ignorespaces}
\def\endproclaim{\rm\ppar}
\def\cite#1{[#1]}
\newcount\itemnumber
\def\roster{\items\itemnumber=1}
\def\endroster{\enditems}
\let\itemold\item
\def\item{\itemold{{\rm(\number\itemnumber)}}%
\global\advance\itemnumber by 1\ignorespaces}
\def\S{section~\ignorespaces}  %%  expand \S to "section"
\def\date#1\enddate{\relax}
\def\thanks#1\endthanks{\relax}   %%%  Move acknowledgements "manually"
\def\dedicatory#1\enddedicatory{\relax}  %%% to the end of the intro.
\def\rom#1{{\rm #1}}  % in some versions of amstex but not all
\let\footnote\plainfootnote

%%
%%
%%   Adapt the amstex reference list macros 
%%   (some spacing may need "manual" adjustment) :
%%
%%
\def\Refs{\ppar{\large\bf References}\ppar\bgroup\leftskip=25pt
\frenchspacing\parskip=3pt plus2pt\small}       
\def\endRefs{\egroup}
\def\widestnumber#1#2{\relax}
\def\endrefitem{}
\def\refdef#1#2#3{\def#1{\leavevmode\unskip\endrefitem#2\def\endrefitem{#3}}}
\def\ref{\par}
\def\endref{\leavevmode\unskip\endrefitem\par\def\endrefitem{}}
\refdef\key{\noindent\llap\bgroup[}{]\ \ \egroup}
\refdef\no{\noindent\llap\bgroup[}{]\ \ \egroup}
\refdef\by{\bf}{\rm, }
\refdef\manyby{\bf}{\rm, }
\refdef\paper{\it}{\rm, }
\refdef\book{\it}{\rm, }
\refdef\jour{}{ }
\refdef\vol{}{ }
\refdef\yr{(}{) }
\refdef\ed{(}{, Editor) }
\refdef\publ{}{ }
\refdef\inbook{from: ``}{'', }
\refdef\pages{}{ }
\refdef\page{}{ }
\refdef\paperinfo{}{ }
\refdef\moreref{}{ }
\refdef\bookinfo{}{ }
\refdef\publaddr{}{ }
\refdef\eds{(}{, Editors) }
\refdef\bysame{\hbox to 3 em{\hrulefill}\thinspace,}{ }
\refdef\toappear{(to appear)}{ }
\refdef\issue{no.\ }{ }
%%
%%      Macros to change refs to numbers
%%
%%      To use these macros uncomment the macros below :
%%
%\newcount\refnumber\refnumber=1
%\def\refkey#1{\expandafter\xdef\csname cite#1\endcsname{\number\refnumber}%
%\global\advance\refnumber by 1}
%\def\cite#1{[\csname cite#1\endcsname]}
%\def\Cite#1{\csname cite#1\endcsname}  %% unbracketed \cite 
%\def\key#1{\noindent\llap{[\csname cite#1\endcsname]\ \ }}
%%
%%      Next edit the reference list so that all keys
%%      are enclosed in { }'s.  Then copy
%%      the list below and edit and cut the lines down to
%%      produce a \refkey list, of which the following 
%%      are sample lines:
%%
%%      \refkey {G2}
%%      \refkey {G3}
%%
%%      Finally check for uses of \cite.  Only simple uses of the key
%%      such as  \cite{G2}  are OK.  More complicated uses such as
%%      \cite{G2; theorem 4.1}  should be edited to use the
%%      unbracketed version \Cite{ },  thus  [\Cite{G2}; theorem 4.1] .
%%
%%       End of patch
%%

\TagsAsMath
%%%%operators:
   \def\card{\#}
   \def\Int{\operatorname{Int}}
   \def\pr{\operatorname{pr}}
%%%%Greek:
   \def\a{\alpha} \redefine\b{\beta} \def\s{\sigma} 
%%%% homemade symbol for connected sums:
\def\twobars#1#2#3#4{\vcenter{\hrule height.#1pt width#2pt
                               \vskip#3pt
                               \hrule height.#1pt width#2pt
                               \vskip#4pt}}
\def\stroke#1#2#3{\vrule height#1pt width.#2pt depth#3pt}
\def\connsum{\hskip2pt
             \twobars4{1.5}31\stroke831\twobars4231\stroke831\twobars4{1.5}31
             \hskip2pt}
%%%%miscellany:
     \def\sub{\subset}
     \redefine\emptyset{\varnothing}
     \def\Bd{\partial}
     \def\C{\Bbb C}
     \def\D{\bold D}
     \def\R{\Bbb R}
     \def\h#1#2{h^{\sssize (#1)}_{#2}} %% handles
     \font\tenss=cmss10 \font\eightss=cmss8
     \def\sanser#1{\mathchoice{\text{\tenss #1}}{\text{\tenss #1}}
                  {\text{\eightss #1}}{\text{\eightss #1}}}
     \def\circ{\sanser O}
     \def\cross{\sanser X}
     \def\o{\sanser o}
     \def\x{\sanser x}
     \def\O#1{\circ_{#1}}
     \def\X#1{\cross_{#1}}
     \font\bxtwelve=cmbx12 \def\PLUMB{\thinspace\botsmash{\lower.75ex
                            \hbox{\text{\bxtwelve*}}}\thinspace}
\hyphenation{Worces-ter dia-gram-mat-ic}
\topmatter
\title
Positive links are strongly quasipositive
\endtitle
\thanks
Partially supported by NSF grant DMS-9504832.
\endthanks
\dedicatory
for Rob Kirby
\enddedicatory
\author
Lee Rudolph
\endauthor
\address
Department of Mathematics, Clark University,
Worcester MA 01610, USA
\endaddress
\abstract
Let $S(\bold D)$ be the surface produced by applying
Seifert's algorithm to the oriented link diagram $\bold D$.
I prove that if $\bold D$ has no negative crossings then
$S(\bold D)$ is a quasipositive Seifert surface, that is,
$S(\bold D)$ embeds incompressibly on a fiber surface plumbed from
positive Hopf annuli. This result, combined with the truth of the
``local Thom Conjecture'', has various interesting consequences;
for instance, it yields an easily-computed estimate for the slice
euler characteristic of the link $L(\bold D)$ (where $\bold D$ is
arbitrary) that extends and often improves the ``slice--Bennequin
inequality'' for closed-braid diagrams; and it leads to yet another
proof of the chirality of positive and almost positive knots.
\endabstract
\asciiabstract{Let S(D) be the surface produced by applying
Seifert's algorithm to the oriented link diagram D.
I prove that if D has no negative crossings then
S(D) is a quasipositive Seifert surface, that is,
S(D) embeds incompressibly on a fiber surface plumbed from
positive Hopf annuli. This result, combined with the truth of the
`local Thom Conjecture', has various interesting consequences;
for instance, it yields an easily-computed estimate for the slice
euler characteristic of the link L(D) (where D is
arbitrary) that extends and often improves the `slice--Bennequin
inequality' for closed-braid diagrams; and it leads to yet another
proof of the chirality of positive and almost positive knots.}
\email
{\tt lrudolph@black.clarku.edu}
\endemail
\primaryclass{57M25}\secondaryclass{32S55, 14H99}
\keywords
Almost positive link,
Murasugi sum,
positive link,
quasipositivity,
Seifert's algorithm
\endkeywords
\makeshorttitle
\cl{\small\it For Rob Kirby}
\document
\head Introduction; statement of results \endhead

Given an oriented link diagram $\D$, let $\X>(\D)$ (resp. $\X<(\D)$)
be the set of positive (resp. negative) crossings,
and $\O\ge(\D)$ (resp. $\O<(\D)$) the set of Seifert circles
adjacent to some $\x\in\X>(\D)$ (resp. to no $\x\in\X>(\D)$).
Write $\card A$ for the number of elements of $A$.
Let $S(\D)$ be the surface produced by applying Seifert's
algorithm to $\D$; let $L(\D):=\Bd S(\D)$.
An oriented link $L$ is
{\it positive} if $L$ is isotopic to $L(\D)$
for some $\D$ with $\card\X<(\D)=0$,
{\it almost positive} if $L$ is not positive but
$L$ is isotopic to $L(\D)$ for some $\D$ with $\card\X<(\D)=1$,
and {\it strongly quasipositive} if $L$ bounds a
{\it quasipositive Seifert surface}, that is, a surface
embedded incompressibly on a fiber surface plumbed
from positive Hopf annuli.

\proclaim{Theorem} If $\card\X<(\D)=0$, then $S(\D)$ is
quasipositive, and so $L(\D)$ is strongly quasipositive.
\endproclaim

Denote by $\chi_s(L)$ the greatest value of the euler
characteristic $\chi(F)$ for $F\sub D^4$ a smooth, oriented
$2$--manifold such that $L=\Bd F$ and $F$ has no closed components.
Let $D^1(K)$ be the untwisted positive Whitehead double of
a knot $K$, $D^k(K):=D^1(D^{k-1}(K))$.

\proclaim{Corollary 1} If $\D$ is any oriented link diagram, then
$$
\chi_s(L(\D))\le %
(\card \O\ge(\D)-\card \O<(\D))-(\card \X>(\D)-\card \X<(\D)).
\tag{*[\D]}
$$
\endproclaim

\proclaim{Corollary 2} \rom{(A)} A non-trivial positive link is chiral.
\rom{(B)} An almost positive knot is chiral.
\endproclaim
This corollary is included chiefly for the novelty of the method;
see Remark~4.

\proclaim{Corollary 3} If $K$ is a non-trivial positive knot,
then $\chi_s(D^k(K))\allowmathbreak=\allowmathbreak-1$, $k>0$.
\endproclaim
This corollary, with ``strongly quasipositive'' in place of
``positive'', was proved in \cite{14}, so it is immediate from
the Theorem.  It is a partial extension of a result of Cochran
and Gompf \cite{1}, which assumes much less than positivity of $K$,
but concludes only that $\chi_s(D^k(K))=-1$, $1\le k\le 6$;
see Remark~5.

\remark{Remarks} (1)\qua As defined above,
strong quasipositivity is an intrinsic geometric property.
Its original definition (\cite{11}; cf \cite{9}), like those
of positivity and almost-positivity, was ``diagrammatic'':
an oriented link $L$ is strongly quasipositive if and only if,
for some $n$, $L$ can be represented as the closure $\widehat\b$
of a braid $\b\in B_n$ that is the product of ``embedded positive bands''
$\s_{i,j}:=(\s_i\dotsc\s_{j-2})\s_{j-1}(\s_i\dotsc\s_{j-2})^{-1}$
(where $\s_1,\dots,\s_{n-1}$ are the standard generators of $B_n$).
The equivalence of this to the intrinsic definition follows by
combining \cite{12} and \cite{16}.

\medskip
{\bf Question}\qua Can positive
links be characterized as strongly quasipositive links that
satisfy some extra geometric conditions?

\medskip
(2)\qua If Corollary~1 is weakened by restricting $\D$ to be a
closed braid diagram and by replacing \thetag{$*[\D]$} by
$$
\chi_s(L(\D))\le %
(\card \O\ge(\D)+\card \O<(\D))-(\card \X>(\D)-\card \X<(\D)),
\tag{*'[\D]}
$$
then it becomes the {\it slice--Bennequin inequality} sBi, \cite{14};
like the proof of Corollary~1 below, the proof of sBi in \cite{14}
makes essential use of the truth of the ``local Thom Conjecture'',
a result originally established using gauge theory for embedded
surfaces \cite{3} which now follows from more general results
established using monopole methods \cite{4}, \cite{5}.
(In fact, \thetag{$*[\D]$} for a closed braid diagram $\D$
already follows from sBi, \cite{16}, Corollary~5.2.2.)

\medskip
(3)\qua A different weakening of Corollary~1---allowing $\D$
to be arbitrary, but concluding only \thetag{$*'[\D]$}---follows
easily by combining sBi and the (proof of the) main
result of \cite{20}.  In fact, after posting the first
version of the present article to the xxx~Mathematics Archives,
I received e-mail from Takuji Nakamura informing me
that Nakamura's February 1998 master's thesis at Keio
University (Japan) gives a proof of the Theorem using
the techniques of \cite{20}; and another reader has since
kindly shown me how to use Vogel's algorithm \cite{19} to
give yet another proof.

\medskip
(4)\qua Chirality of non-trivial positive links was
established well over a decade ago: the proof in \cite{1}
(which does not, itself, depend on ``Donaldson's theorems'',
ie, gauge theory, but uses only a classical invariant,
the signature) was described by Cochran and Gompf in talks
at MSRI in 1985 and at the 1986 Santa Cruz Summer Research
Conference on braids; other proofs using the signature
were published independently by Przytycki \cite{7}
and Traczyk \cite{18} (cf \cite{10} for an earlier special case,
with an unconventional choice of sign).
Chirality of almost positive links (via negativity of the signature)
was announced in a 1991 abstract \cite{8} by Przytycki and Taniyama;
Stoimenow has recently given a proof \cite{17} (for knots) which uses
a Vassiliev-style invariant due to Fiedler.
Chirality of non-trivial strongly quasipositive knots
follows from a calculation with the FLYPMOTH link polynomial
(\cite{6}, remark following Problem~9), and chirality of
non-trivial strongly quasipositive links is an easy
consequence of sBi (cf \cite{6}, Problems 8.2 and 9.2).
When I came across \cite{17} while
perusing {\tt http://front.math.ucdavis.edu/math.GT/} preparatory
to uploading what I had thought was the finished version of the present
paper, I realized that an alternative proof of the chirality of
almost positive links could be given using Corollary~1,
and I revised this paper accordingly.
Many thanks are due to the maintainers both of the xxx~Mathematics
Archives and of the Front for keeping us all on our toes.
Thanks also to Tim Cochran, Jozef Przytycki, and others
for helpful e-mail on the history of chirality results for positive
and almost positive links.

\medskip
(5)\qua A long-standing conjecture (\cite{2}, Problem 1.38)
asserts that, for a knot $K$, $\chi_s(D^1(K))=1$ iff $\chi_s(K)=1$.
Since any Whitehead double of a knot bounds a
punctured torus (already in $S^3$), this conjecture 
implies a second conjecture: if $\chi_s(K)<1$, then
$\chi_s(D^k(K))=-1$ for $1\le k < \infty$.
Cochran and Gompf \cite{1} made some progress towards
the second conjecture: they defined what it means for a knot
$K$ to be ``greater than or equal to the positive trefoil''
(briefly, $K\ge\Bbb T$); they proved (by way of providing
a large class of examples) that if $K\ne O$ is a positive knot
(in particular, a closed positive braid) then $K\ge\Bbb T$;
and they applied gauge theory to show that if $K\ge\Bbb T$
then $\chi_s(K)<1$ and $\chi_s(D^k(K))=-1$ for $1\le k\le 6$.
In \cite{14}, I proved the second conjecture if $K\ne O$
is a strongly quasipositive knot (in particular, a closed
positive braid).  Corollary~3 is the observation that the
second conjecture is true for positive knots.

\medskip
{\bf Questions}\qua Is the second conjecture
true for all $K\ge\Bbb T$?  Might it
in fact be the case that $K\ge\Bbb T$ implies $K$ is strongly
quasipositive?

\medskip
(6)\qua If $L(\D)$ is a knot $K$, then the estimate in \cite{15}
can be rewritten as
$$
\chi_s(K)\le m(\D)-(\card \X>(\D)-\card \X<(\D)-s_-(\D)),
\tag{*''[\D]}
$$
where $\D\sub\C$ is taken to be in general position
with respect to $\Im\co \C\to\R:z\mapsto(z-\bar z)/2i$,
$m(\D)$ is the number of local maxima of $\Im|\D$,
and $s_-(\D)$ is the number of crossings of $\D$
which (disregarding their actual signs) are ``locally
negative'' when oriented by $\Im$.  This inequality
resembles \thetag{$*'[\D]$} (or sBi).

\medskip
{\bf Question}\qua Is there---in general, or in the
special case that $s_-(\D)=0$ (so that $\D$ is a
{\it positive plat diagram} of $K$, \cite{13})---a modification
of \thetag{$*''[\D]$}, analogous to \thetag{$*[\D]$}, in which
(for some $\D$) some local maxima of $\Im|\D$ contribute $-1$
rather than $1$ to the right-hand side?

\medskip
(7)\qua Let $u(K)$ denote the {\it unknotting number} of the knot $K$.
Of course for all $K$, $u(K)\ge (1-\chi_s(K))/2$; there exist $K$
such that the estimate for $u(K)$ based on \thetag{$*[\D]$}
is sharper than the estimate based on sBi.
\endremark
\medskip

\head Proof that positive links are strongly quasipositive\endhead

In preparation for the proofs, recall Seifert's algorithm.
Given an oriented link diagram $\D\sub\C$
(where $\C$ has its standard orientation),
let $\cross(\D):=\X>(\D)\cup\X<(\D)$,
$\circ(\D):=\O\ge(\D)\cup\O<(\D)$.
The algorithm, given input $\D$, produces
a piecewise-smooth oriented surface $S(\D)\sub\C\times\R$
equipped with a $(0,1)$--handle decomposition
$S(\D)=\bigcup_{\o\in \circ(\D)}\h0\o\cup \bigcup_{\x\in\cross(\D)}\h1\x$
such that:

\roster{\leftskip 5em
\itemold{(1)~\qquad} \hskip-2.5em
$\pr_1(\Bd S(\D))$ is the underlying
oriented immersed $1$--manifold of $\D$;
\itemold{(2)~\qquad} \hskip-2.5em
for every $\o\in\circ(\D)$, $\pr_1|\h0\o\co \h0\o\to\C$ is a
(non-oriented) embedding with $\pr_1(\Bd\h0\o)=\o$;
\itemold{(3)~\qquad} \hskip-2.5em
for every $\x\in\cross(\D)$,
\itemold{(3.1)}
there is one transverse arc
$\tau(\h1\x)\sub\h1\x$ such that $\pr_1(\tau(\h1\x))=\x$,
\itemold{(3.2)}
$\pr_1|(\h1\x\setminus\tau(\h1\x))$ is an embedding,
and preserves orientation on precisely one component of
$\h1\x\setminus\tau(\h1\x)$,
\itemold{(3.3)}
$\h1\x$ is attached to $\h0\o$ iff $\x$ is adjacent to $\o$ in $\D$, and
\itemold{(3.4)}
the over-arc of $\D$ through $\x$ contains
the image by $\pr_1$ of the component of $\h1\x\cap\Bd S(\D)$
which contains that endpoint of $\tau(\h1\x)$ at which $\pr_2$
takes on the larger value.\par}
\endroster
Call $\o\in\circ(\D)$ {\it outermost} if
$\o\cap\pr_1(\h0{\o'})=\emptyset$ for $\o\ne\o'\in\circ(\D)$;
let $\circ\circ(\D)$ be the set of outermost Seifert circles
of $\D$.  Of course $\card\circ\circ(\D)\ge 1$.

\demo{Proof of Theorem} Let $\card\X<(\D)=0$.
The split sum of quasipositive Seifert surfaces is
quasipositive, so there is no loss of generality in
assuming that $\D$ is connected.  The proof proceeds
by induction, first on $\card\cross(\D)$ and then
on $|\card\circ\circ(\D)-2|$.

If $\card\cross(\D)=0$ then $\D$ is the trivial diagram
of the trivial knot and the Theorem is trivially true.

If $\card\cross(\D)>0$ and $\card\circ\circ(\D)=1$,
then $\D$ can be replaced by $\D'$ (isotopic to $\D$
on $S^2=\C\cup\{\infty\}$) such that
$\circ(\D')=\circ(\D)\cup\{\o'_0\}\setminus\{\o_0\}$,
$\cross(\D')=\cross(\D)$, and $S(\D')$ is ambient isotopic to
$S(\D)$ in $\C\times\R$.  Thus it may be assumed
that $\card\circ\circ\D\ge 2$.  Since $\D$ is connected,
there exist $\o_1, \o_2\in\circ\circ(\D))$ such that
$\pr_1|\h0\o\co \h0{\o_1}\to\C$ preserves orientation and
$\pr_1|\h0\o\co \h0{\o_2}\to\C$ reverses orientation.

If $\card\circ\circ(\D)=2$, then the
union $F$ of $\h0{\o_1}$ and $\h0{\o_2}$, together with
all the $1$--handles $\h1\x$ such that $\x$ is adjacent to
both $\o_1$ and $\o_2$, is ambient isotopic to $S(\D(\s_1^k))$,
where $k\ge 1$ is the number of those $1$--handles
and $\D(\s_1^k)$ is the positive closed braid diagram
of $\s_1^k\in B_2$; moreover, for $i=1,2$, the union
$G_i$ of $\h0{\o_i}$ with all the $0$--handles $\h0\o$
such that $\o\sub\pr_1(\h0{\o_i})$, together with all
the $1$--handles $\h1\x$ to which any of these $0$--handles
is attached, is ambient isotopic to $S(\D_i)$ for an appropriate
positive oriented link diagram; finally, $S(\D)$ is an iterated
Murasugi sum (Stallings plumbing) $G_1\PLUMB F \PLUMB G_2$.
Now, $S(\D(\s_1^k))$, the fiber surface of the $(2,k)$ torus link,
is well known to be quasipositive (an explicit plumbing
from positive Hopf annuli can be extracted from \cite{12}),
while $S(\D_1)$ and $S(\D_2)$ are quasipositive by induction
on the number of crossings.  Since, by \cite{16},
any plumbing of quasipositive Seifert surfaces is quasipositive,
$S(\D)$ is quasipositive.

If $\card\circ\circ(\D)>2$, then there exist
$\o', \o''\in\circ\circ(\D)$, $\o'\ne\o''$, such that
$\pr_1|\h0{\o'}$ and $\pr_1|\h0{\o''}$ have the same
orientation type (both preserve, or both reverse, orientation);
a little thought about planar embeddings of bipartite graphs
shows that such $\o', \o''$ may be chosen to have the further property
that there is an arc $\a\sub\C$ with one endpoint on $\o'$ and the
other endpoint on $\o''$, which is otherwise disjoint from $\D$.
Let $\D_\a$ be the positive link diagram
with $\cross(\D_\a)=\cross(\D)$ and
$\circ(\D_\a)=\circ(\D)\cup\{\o'\connsum_\a\,\o''\}\setminus\{\o',\o''\}$,
where
$\o'\connsum_\a\,\o'':=\Bd(\pr_1(\h0{\o'})\cup N(\a)\cup\pr_1(\h0{\o''}))$
for a suitable relative regular neighborhood $N(\a)$ of $\a$
in $\C\setminus\Int(\pr_1(\h0{\o'}\cup\h0{\o''}))$.
Evidently $S(\D_\a)$ may be constructed to contain $S(\D)$,
and then $S(\D)$ is clearly embedded incompressibly on $S(\D_\a)$.
By induction on the number of outermost Seifert circles,
$S(\D_\a)$ is quasipositive.  Since, by \cite{12},
an incompressible subsurface of a quasipositive Seifert surface
is quasipositive, it follows that $S(\D)$ is quasipositive.
\qed\enddemo

\head Proofs of the corollaries\endhead

\demo{Proof of Corollary 1} Let
$Q(\D):=\bigcup_{\o\in \O\ge(\D)}\h0\o\cup\bigcup_{\x\in\X>(\D)}\h1\x
\sub S(\D)$.
Then $Q(\D)$ is clearly $S(\D^{+})$, where
$\pr_1(\Bd Q(\D))$ is the underlying oriented immersed
$1$--manifold of $\D^{+}$, $\circ(\D^+)=\O\ge(\D^+)=\O\ge(\D)$, and
$\cross(\D^+)=\X>(\D^+)=\X>(\D)$.  By the Theorem, $Q(\D)$ is
quasipositive, so
$$
\align
\chi_s(L(\D))
         &\le 2\chi(Q(\D))-\chi(S(\D)) \\
         &=2(\card \O\ge(\D)-\card \X>(\D))-(\card\circ(\D)-\card\cross(\D))\\
         &=(\card \O\ge(\D)-\card \O<(\D))-(\card \X>(\D)-\card \X<(\D))
\endalign
$$
by \cite{16}, Corollary~5.2.2.
\qed\enddemo

\demo{Proof of Corollary 2}
(A)\qua As noted in Remark~3, chirality of non-trivial strongly quasipositive
links is a consequence of sBi, so by the Theorem, if
$L$ is a non-trivial positive link then $L$ is chiral.

\medskip
(B)\qua Let $\D$ be an oriented link diagram such that $L(\D)=K$ is a knot
and $\card\X<(\D)=1$.  Then $\card\O<(\D)\le 1$.  The following
case analysis shows that either $K$ is positive (so that
$K$ is not almost positive) or $\chi_s(K\connsum K)\le -1$ (so that
$K$ is not amphicheiral). 

(1)\qua If $\card\O<(\D)=1$ then a single
Reidemeister move of type $1$ replaces $\D$ by $\D_0$ with
$L(\D_0)=K$ and $\card\X<(\D)=0$, so $K$ is positive.

(2)\qua If $\card\O<(\D)=0$ and $\chi(S(\D))=\card\O>(\D)-(\card\X>+1)<-1$,
then Corollary~1, applied to the connected sum $\D\connsum\D$,
implies that $\chi_s(K\connsum K)\le -1$.

(3)\qua If $\chi(S(\D))=-1$, then a simple
analysis of the possible configurations of the Seifert circles
shows (without using any information about the signs of the
crossings) that the punctured torus $S(\D)$ is isotopic
to either 

(i)\qua a plumbing $A(O,p)\PLUMB A(O,q)$ of two unknotted
twisted annuli, or 

(ii)\qua a pretzel surface $P(2a+1,2b+1,2c+1)$;
if, further, $\card\X<(\D)=1$, then in case~(i) $p$ and
$q$ are non-positive so $K$ is positive by inspection, while
in case~(ii) two of $2a+1,2b+1,2c+1$ are negative and the third
is $1$, whence $P(2a+1,2b+1,2c+1)$ is again isotopic to a plumbing
$A(O,p)\PLUMB A(O,q)$ with $p$ and $q$ non-positive
(if $a,b<-1$, $c=1$, then $p=a+1$, $q=b+1$), so again $K$ is positive.
\qed\enddemo

\rk{Acknowledgement} The author is partially supported by 
NSF grant DMS-9504832.

\Addresses\recd

\bye